\documentclass[12pt]{article}
\usepackage{latexsym,amsmath,amsthm}
\usepackage{amssymb}
\usepackage{graphicx}
\usepackage{appendix}
\usepackage[dvipsnames]{xcolor}
\usepackage{tikz}
\usepackage{setspace}

\newtheorem{thm}{Theorem}[section]

\newtheorem{conj}[thm]{Conjecture}
\newtheorem{claim}[thm]{Claim}

\theoremstyle{definition}

\theoremstyle{definition}

\theoremstyle{definition}

\theoremstyle{definition}

\theoremstyle{definition}

\theoremstyle{definition}

\theoremstyle{remark}

\theoremstyle{remark}

\usepackage[margin=1in]{geometry}

\usepackage{url}
\usepackage{microtype}

\begin{document}

\title{The domination polynomial of powers of paths and cycles}

\author{David Galvin\thanks{Department of Mathematics,
University of Notre Dame, Notre Dame IN USA; {\tt dgalvin1@nd.edu}. Galvin in part supported by a Simons Collaboration Grant for Mathematicians.} \and Yufei Zhang\thanks{Department of Mathematics,
University of Notre Dame, Notre Dame IN USA.}}

\maketitle

\begin{abstract}
A dominating set in a graph is a set of vertices with the property that every vertex in the graph is either in the set or adjacent to something in the set. The domination sequence of the graph is the sequence whose $k$th term is the number of dominating sets of size $k$.

Alikhani and Peng have conjectured that the domination sequence of every graph is unimodal. Beaton and Brown verified this conjecture for paths and cycles. Here we extend this to arbitrary powers of paths and cycles.  
\end{abstract}

\section{Introduction}

A {\em dominating set} in a graph $G$ is a subset $D$ of vertices of $G$ with the property that every vertex in $G$ is either in $D$ or adjacent to something in $D$. Write $\gamma_k(G)$ for the number of dominating sets of size $k$ in $G$. For a thorough survey of dominating sets, see the monograph by 
Haynes,  Hedetniemi and  Slater \cite{HaynesHedetniemiSlater1998}.

In \cite{AlikhaniPeng2014} Alikhani and Peng propose the following conjecture. Here a sequence $(a_k)_{k\geq 0}$ is {\it unimodal} if there is an index $m$ --- a {\it mode} of the sequence --- with 
$$
a_0 \leq a_1 \leq \cdots \leq a_m \geq a_{m+1} \geq a_{m+2} \cdots.
$$
Note that a unimodal sequence does not necessarily have a unique mode.
\begin{conj} \label{conj-dom-unimodal}
For all graphs $G$ the sequence $(\gamma_k(G))_{k\geq 0}$ is unimodal.
\end{conj}
There are partial results in the direction of Conjecture \ref{conj-dom-unimodal} (see below) but it remains open in general. As observed by Beaton and Brown \cite{BeatonBrown2022} the domination sequence of a graph does not in general satisfy the stronger property of log-concavity (a sequence $(a_k)_{k \geq 0}$ is {\em log-concave} if $a_k^2 \geq a_{k-1}a_{k+1}$ for all $k \geq 1$). This may partly explain why Conjecture \ref{conj-dom-unimodal} has proven tricky; there are fewer techniques available to establish unimodality than log-concavity.   

We now briefly mention some of the recent partial results on unimodality of the domination sequence (for some less recent result see Beaton's survey \cite{Beaton2017}). Beaton and Brown \cite{BeatonBrown2022} settle the conjecture in the affirmative for complete multipartite graphs (and other families --- see below). Burcroff and O'Brien \cite{BurcroffOBrien2023} establish the conjecture for spiders (trees with at most one vertex of degree greater than two) whose maximum degree is at most 400; for lollipop graphs (complete graphs with a path dropped from one vertex), direct products of complete graphs, and Cartesian products of two complete graphs. S. Zhang \cite{Zhang2021} shows that if $G$ has a universal vertex (one adjacent to all other vertices) and more than $2^{13}$ vertices then its domination sequence is unimodal. 

Alikhani and Peng \cite{AlikhaniPeng2010} show that the domination sequence of every $n$-vertex graph is weakly increasing up to $\lceil n/2 \rceil$. In \cite{Beaton2021} Beaton shows that if an $n$-vertex graph has no isolated vertices then its domination sequence is weakly decreasing from $\lfloor 3n/4 \rfloor$ on. Burcroff and O'Brien \cite{BurcroffOBrien2023} show that if $n$-vertex $G$ has $m$ universal vertices then its domination sequence is weakly decreasing from $\lceil n/2 + n/2^{m+1} \rceil$ on.   

Two results --- Theorems \ref{thm-from-bb-min-degree} and \ref{thm-BB-path-cycle} below, both due to Beaton and Brown \cite{BeatonBrown2022} --- are particularly relevant to the present paper. 
\begin{thm} \label{thm-from-bb-min-degree}
If $G$ is a graph with $n$ vertices and with minimum degree at least $2\log_2n$ then the domination sequence of $G$ is unimodal with mode at $\lceil n/2 \rceil$.
\end{thm}
From Theorem \ref{thm-from-bb-min-degree} Beaton and Brown deduce that for any $p \in (0,1)$ the probability that the Erd\H{o}s-R\'enyi random graph $G(n,p)$ has unimodal domination sequence tends to $1$ as $n \rightarrow \infty$ (i.e., ``Conjecture \ref{conj-dom-unimodal} is almost always true'').  
\begin{thm} \label{thm-BB-path-cycle}
For every $n \geq 0$, both the path on $n$ vertices and the cycle on $n$ vertices have unimodal domination sequences. 
\end{thm}

The main point of this note is to extend Theorem \ref{thm-BB-path-cycle}. For a graph $G$ and a positive integer $\ell$ the {\em $\ell$th power} of $G$ is the graph on the same vertex set as $G$, with two vertices adjacent if their distance (in $G$) is at most $\ell$ (so the first power of a graph is the graph itself). Denote by $P^\ell_n$ the $\ell$th power of the path on $n$ vertices, and by $C^\ell_n$ the $\ell$th power of the cycle on $n$ vertices. We interpret both $P^\ell_0$ and $C^\ell_0$ to be the graph with no vertices (this graph has one dominating set, of size $0$). We interpret $C^\ell_1$ to be a single isolated vertex and $C^\ell_2$ to the complete graph on two vertices. 

\begin{thm} \label{thm-powers-of-path-cycle}
For every $n \geq 0$ and $\ell \geq 1$, both $P^\ell_n$ and $C^\ell_n$ have unimodal domination sequences. 
\end{thm}

Denote by $\gamma(G,x)$ the {\em domination polynomial} of graph $G$, that is, the polynomial 
$$
\gamma(G,x) = \sum_{k \geq 0} \gamma_k(G)x^k.
$$
One step in the proof of Theorem \ref{thm-powers-of-path-cycle} involves recurrence relations for $\gamma(P^\ell_n,x)$ and $\gamma(C^\ell_n,x)$.
\begin{thm} \label{thm-power-recurrences}
For $\ell \geq 1$ and $n \geq 2\ell+1$
\begin{equation} \label{path-power-recurrence}
\gamma(P^\ell_n,x) = x\sum_{j=1}^{2\ell+1} \gamma(P^\ell_{n-j},x).
\end{equation}
For $\ell \geq 1$ and $n \geq 2\ell+2$
\begin{equation} \label{cycle-power-recurrence}
\gamma(C^\ell_n,x) = x\sum_{j=1}^{2\ell+1} \gamma(C^\ell_{n-j},x).
\end{equation}
\end{thm}
The case $\ell=1$ of \eqref{path-power-recurrence} was first established by Alikhani and Peng in \cite{AlikhaniPeng2009}, via a quite involved argument, and a much simplified proof was given by Arocha and Llano in \cite{ArochaLlano2016}. The cases $\ell=2, 3$ of \eqref{path-power-recurrence} were obtained by Vijayan and Gibson \cite{VijayanGibson2013} and Medone and Christilda \cite{MedoneChristilda2022}, also by quite long arguments. The case $\ell=1$ of \eqref{cycle-power-recurrence} appears in \cite{BeatonBrown2022}. Our proofs of \eqref{path-power-recurrence} and \eqref{cycle-power-recurrence} for general $\ell$ are quite short and direct.

The recurrence relations \eqref{path-power-recurrence} and \eqref{cycle-power-recurrence} will take center stage in the proof of Theorem \ref{thm-powers-of-path-cycle}. The proof will also require us to understand the initial conditions of these recurrences, so it will be useful to have a way to easily compute $\gamma(P^\ell_n,x)$ and $\gamma(C^\ell_n,x)$ for values of $n$ that are small compared to $\ell$. For cycles, this is very easy. For $n \leq 2\ell+1$ we have that $C^\ell_n$ is the complete graph on $n$ vertices, for which every non-empty subset of vertices is dominating, and so
\begin{equation} \label{eq-small-n-cycle-recurrence}
\gamma(C^\ell_n,x) = (1+x)^n-1.
\end{equation}
For paths things are not as simple, since for $\ell +2 \leq n \leq 2\ell$ the graph $P^\ell_n$ is not a complete graph. Nevertheless we have two ways to compute $\gamma(P_n^\ell,x)$ for small values of $n$. The first is a direct extension of \eqref{path-power-recurrence}, while the second uses \eqref{path-power-recurrence} as a starting point but then goes in a different direction.

\begin{thm} \label{thm-small-n-power-recurrences}
For $n \geq 0$ and $\ell \geq 1$ we can compute $\gamma(P_n^\ell,x)$ via the following recurrence:
\begin{description}
\item[A1] $\gamma(P_0^\ell,x)=1$ for $\ell \geq 1$.
\item[A2] $\gamma(P_1^\ell,x)=x$ for $\ell \geq 1$.
\item[A3] $\gamma(P_n^\ell,x) = nx + x\sum_{j=1}^{n-1} \gamma(P_{n-j}^\ell,x)$ for $2 \leq n \leq \ell+1$.
\item[A4] $\gamma(P_n^\ell,x) = (2\ell+2-n)x + x\sum_{j=1}^{n-1} \gamma(P_{n-j}^\ell,x)$ for $\ell+2 \leq n \leq 2\ell$.
\item[A5] $\gamma(P_n^\ell,x) = x\sum_{j=1}^{2\ell+1} \gamma(P_{n-j}^\ell,x)$ for $2\ell+1 \leq n$ (note that this is just \eqref{path-power-recurrence}).
\end{description}

Alternatively for $n \geq 0$ and $\ell \geq 1$ we can compute $\gamma(P_n^\ell,x)$ via
\begin{description}
\item[B1] $\gamma(P_0^{\ell},x)=1$ for $\ell \geq 1$.
\item[B2] $\gamma(P_n^{\ell},x) = (1+x)^n -1$ for $1 \leq n \leq \ell+1$.
\item[B3] $\gamma(P_n^{\ell},x) = (1+x)\gamma(P_{n-1}^{\ell},x) - x$ for $\ell+2 \leq n \leq 2\ell+1$.
\item[B4] $\gamma(P_n^{\ell},x) = (1+x)\gamma(P_{n-1}^{\ell},x) - x\gamma(P_{n-2(\ell+1)}^{\ell},x)$ for $n \geq 2\ell+2$. 
\end{description}
\end{thm}

A second step in the proof of Theorem \ref{thm-powers-of-path-cycle} involves a generalization of a result of Beaton and Brown \cite{BeatonBrown2022}. Say that a polynomial is {\it unimodal} if its coefficient sequence is unimodal, and say that a sequence $(a_k)_{k=c}^d$ is {\it barely increasing} if $0 \leq a_{k+1}-a_k \leq 1$ for all $c \leq k \leq d-1$.
\begin{thm} \label{thm-bb-gen}
Fix $k \geq 3$. Let $f_0(x), f_1(x), \ldots$ be a sequence of polynomials satisfying 
$$
f_n(x) = x \sum_{j=1}^k f_{n-j}(x)
$$
for all $n \geq k$. Let ${\mathcal P}_n$ be the property that 
\begin{itemize}
\item the coefficient sequence of $f_i$ is non-negative and unimodal for each $0 \leq i \leq n$, and
\item the $f_i$'s have a sequence of modes $(m_i)_{i=0}^n$ that is barely increasing. 
\end{itemize}
If ${\mathcal P}_k$ holds then ${\mathcal P}_n$ holds for all $n \geq k$.   
\end{thm}
The case $k=3$ of Theorem \ref{thm-bb-gen} appears in \cite{BeatonBrown2022}. Note that we cannot replace ${\mathcal P}_k$ with ${\mathcal P}_{k-1}$ in Theorem \ref{thm-bb-gen}. For example, setting $f_0(x)=3$, $f_1(x)=x$ and $f_i(x)=2x^2$ for $2 \leq i \leq k-1$, we have that ${\mathcal P}_{k-1}$ holds but not ${\mathcal P}_k$.

The final ingredient needed for the proof of Theorem \ref{thm-powers-of-path-cycle} is Theorem \ref{thm-from-bb-min-degree}, introduced earlier.

We briefly describe the layout of the paper here.
\begin{itemize}
\item We establish the first part of Theorem \ref{thm-power-recurrences} (\eqref{path-power-recurrence}, the recurrence for the domination polynomial of powers of paths) in Section \ref{subsec-proofs-power-recurrences-path-1}. The proof of the second part of Theorem \ref{thm-power-recurrences} (\eqref{cycle-power-recurrence}, the recurrence for the domination polynomial of powers of cycles) is given in Section \ref{subsec-proofs-power-recurrences-cycle}. 
\item The proof of Theorem \ref{thm-small-n-power-recurrences} (the recurrences for the domination polynomials of powers of paths when $n$ is small compared to $\ell$) is presented in Section \ref{subsec-proofs-power-recurrences-path-3}.
\item In Section \ref{subsec-proofs-power-recurrences-path-2} we make a digression to give a cleaner verification of a slightly weaker version of the path recurrence \eqref{path-power-recurrence}.
\item The proof of Theorem \ref{thm-bb-gen} (on unimodality of sequences of polynomials) is given in Section \ref{subsec-proofs-bb-gen}. 
\item The proof of the main result, Theorem \ref{thm-powers-of-path-cycle} (unimodality of domination sequence for powers of paths and cycles) appears in Section \ref{subsec-proofs-powers-of-path-cycle}.
\end{itemize}

To conclude the introduction, we note that while the domination sequence of a graph is not in general log-concave, it may be for powers of paths and cycles. Indeed, we have checked via {\tt Mathematica} that both $(\gamma(P^\ell_n,k))_{k \geq 0}$ and $(\gamma(C^\ell_n,k))_{k \geq 0}$ are log-concave for all $n, \ell \leq 500$, and even satisfy the stronger property of ultra log-concavity. (A sequence $(a_k)_{k=0}^n$ is {\it ultra log-concave} if $(a_k/\binom{n}{k})_{k=0}^n$ is log-concave.)
\begin{conj} \label{conj-log-concave}
For all $n \geq 0$ and $ \ell\geq 1$ both of $(\gamma(P^\ell_n,k))_{k \geq 0}$ and $(\gamma(C^\ell_n,k))_{k \geq 0}$ are ultra log-concave.
\end{conj}

\section{Proofs} \label{sec-proofs}

\subsection{Recurrences for powers of paths and cycles} \label{subsec-proofs-power-recurrences} 

\subsubsection{The main recurrence for powers of paths} \label{subsec-proofs-power-recurrences-path-1} 

We prove \eqref{path-power-recurrence} combinatorially, starting by considering $n \geq 2\ell+2$.

Let the path $P_n$ underlying $P_n^\ell$ have vertices $v_1, v_2, \ldots, v_n$, with $v_i$ adjacent to $v_{i+1}$ for $i=1, \ldots, n-1$. Put an order on the vertices in the natural way, by saying that $v_i$ is smaller than $v_j$ iff $i < j$.

A dominating set $D$ in $P^\ell_n$ must have at least two vertices. To ensure that $v_1$ is either in $D$ or adjacent to something in $D$, the smallest vertex of $D$ must be one of $\{v_1, \ldots, v_{\ell+1}\}$, say $v_f$. To ensure that $v_{f+\ell+1}$ is either in $D$ or adjacent to something in $D$, the second smallest vertex of $D$ must be one of $\{v_{f+1}, \ldots, v_{f+2\ell+1}\}$, say $v_{f+s}$. For each $f=1, \ldots, \ell+1$ and $s=1, \ldots, 2\ell+1$, denote by ${\mathcal P}_n^\ell(f,s)$ the set of dominating sets of $P_n^\ell$ in which the first vertex is $v_f$ and the second is $v_{f+s}$. Note that ${\mathcal P}_n^\ell(f,s)$ will be empty if $f+s > n$ and non-empty otherwise.

For each $s=1, \ldots, 2\ell+1$ let the path $P_{n-s}$ underlying $P_{n-s}^\ell$ have vertices $v^s_1, v^s_2, \ldots, v^s_{n-s}$, with $v^s_i$ adjacent to $v^s_{i+1}$ for $i=1, \ldots, n-s-1$, and put the natural linear order on the vertices as before. To ensure that $v^s_1$ is either in the dominating set or adjacent to something in the dominating set, a dominating set in $P^\ell_{n-s}$ must have at least one vertex from among $\{v^s_1, \ldots, v^s_{\ell+1}\}$. Denote by ${\mathcal P}_{n-s}^\ell(f)$ the set of dominating sets of $P_{n-s}^\ell$ in which the first vertex is $v^s_f$. Note that ${\mathcal P}_{n-s}^\ell(f)$ will be empty if $f+s > n$ and non-empty otherwise. 

\begin{claim} \label{clm-path-bijection}
For $f+s \leq n$ (when both ${\mathcal P}_n^\ell(f,s)$ and ${\mathcal P}_{n-s}^\ell(f)$ are non-empty) there is a bijection $m_n^\ell(f,s):{\mathcal P}_n^\ell(f,s) \rightarrow {\mathcal P}_{n-s}^\ell(f)$ that reduces sizes of sets by $1$. 
\end{claim}

Claim \ref{clm-path-bijection} would imply that
\begin{eqnarray*}
\gamma(P_n^\ell,x) & = & [f+s \leq n]\sum_{s=1}^{2\ell+1} \sum_{f=1}^{\ell+1} \sum_{D \in {\mathcal P}_n^\ell(f,s)} x^{|D|} \\
& = & [f+s \leq n]\sum_{s=1}^{2\ell+1} x \sum_{f=1}^{\ell+1} \sum_{D \in {\mathcal P}_{n-s}^\ell(f)} x^{|D|} \\
& = & x\sum_{s=1}^{2\ell+1} \gamma(P_{n-s}^\ell,x)
\end{eqnarray*}
which is \eqref{path-power-recurrence}. Here $[\cdot]$ is the Iverson bracket that takes the value $1$ if $\cdot$ is true, and takes the value $0$ otherwise. 

\begin{proof} (Proof of Claim \ref{clm-path-bijection})
Define a map $m_n^\ell(f,s)$ from ${\mathcal P}_n^\ell(f,s)$ to subsets of the vertex set of $P_{n-s}^\ell$ as follows. Let $D=\{v_f, v_s, v_{s_1}, v_{s_2}, \ldots, v_{s_k}\}$ (with $s < s_1 < s_2 < \cdots < s_k$) be a dominating set of $P_n^\ell$ of size $k+2$ ($k \geq 0$) that lies in ${\mathcal P}_n^\ell(f,s)$. Set
$$
m_n^\ell(f,s)(D)=\{v^s_f, v^s_{s_1-(s-f)}, v^s_{s_2-(s-f)}, \ldots, v^s_{s_k-(s-f)}\}.
$$ 
Informally this map deletes the vertices between $v_f$ and $v_s$ and identifies $v_f$ and $v_s$, to go from $P_n^\ell$ to $P_{n-s}^\ell$. It leaves $D$ unchanged, except that the identification of $v_f$ and $v_s$ drops the size of $D$ by $1$. 

Evidently $\left|m_n^\ell(f,s)(D)\right|=|D|-1$. Evidently also $m_n^\ell(f,s)$ is an injective map (if $D \neq D' \in {\mathcal P}_n^\ell(f,s)$, and $v_t$ is in one of $D, D'$ but not in the other, then $v^s_{t-(s-f)}$ is in one of $m_n^\ell(f,s)(D)$, $m_n^\ell(f,s)(D')$ but not the other). 

We claim that $m_n^\ell(f,s)(D)$ is in ${\mathcal P}_{n-s}^\ell(f)$. Clearly $v^s_f$ is the smallest vertex in $m_n^\ell(f,s)(D)$. We now must show that $m_n^\ell(f,s)(D)$ is a dominating set in $P_{n-s}^\ell$. Each of $v^s_1, \ldots, v^s_{f-1}$  is adjacent to $v^s_f$. Now consider a vertex $v^s_t$ in $P_{n-s}^\ell$, for $t > f$. There are three possibilities:
\begin{enumerate}
\item $v_{t+(s-f)}$ (in $P_n^\ell$) is in $D$. In this case $v^s_t \in m_n^\ell(f,s)(D)$.
\item $v_{t+(s-f)} \not \in D$, but $v_{t+(s-f)}$ is adjacent to $v_s$. In this case $v^s_t$ is adjacent to $v^s_f$ (in $P_{n-s}^\ell$).
\item $v_{t+(s-f)} \not \in D$, $v_{t+(s-f)}$ is not adjacent to $v_s$, but $v_{t+(s-f)}$ is adjacent to some $v_{s'} \in D$. In this case $v^s_t$ is adjacent to $v^s_{s'-(s-f)}$ in $P_{n-s}^\ell$, and $v^s_{s'-(s-f)} \in m_n^\ell(f,s)(D)$.  
\end{enumerate}
This shows that indeed $m_n^\ell(f,s)(D)$ is a dominating set in $P_{n-s}^\ell$ and so is in ${\mathcal P}_{n-s}^\ell(f)$. 

To complete the proof of Claim \ref{clm-path-bijection} we need to show that $m_n^\ell(f,s)(D)$ is surjective. Let $D'=\{v^s_f, v^s_{s_1}, \ldots, v^s_{s_k}\}$ (with $f < s_1 < \cdots < s_k$) be in ${\mathcal P}_{n-s}^\ell(f)$. Set 
$$
D = \{v_f, v_s, v_{s_1+(s-f)}, \ldots, v_{s_k-(s-f)}\}.
$$ 
By almost identical reasoning to the argument in the last paragraph we have that $D \in {\mathcal P}_n^\ell(f,s)$; the one point that needs to be added is that each of $v_{f+1}, \ldots, v_{s-1}$ in $P_n^\ell$ is adjacent to at least one of $v_f, v_s$, by the condition that $s \leq 2\ell+1$. The proof is completed by noting that it is evident that $m_n^\ell(f,s)(D)=D'$. 
\end{proof}

The combinatorial proof just presented is only valid when $n \geq 2\ell+2$. Modifying it slightly we can extend the recurrence to be valid for smaller values of $n$. We start with $n=2\ell+1$ (thus completing the proof of \eqref{path-power-recurrence}). Note that for this value of $n$ it is no longer the case that every dominating set of $P_n^\ell$ has at least two elements. Specifically, there is exactly one dominating set, namely $\{v_{\ell+1}\}$, that has size $1$, and this set contributes $x$ to $\gamma(P_{2\ell+1}^\ell,x)$.

As before, for each $s=1, \ldots, 2\ell$ we have that 
$$
\bigcup_{f:(f,s) \neq \emptyset} m_{2\ell+1}^\ell(f,s)\left({\mathcal P}_{2\ell+1}^\ell(f,s)\right)
$$
is the set of all dominating sets of $P^\ell_{2\ell+1-s}$. However, this is no longer the case when $s=2\ell+1$, since ${\mathcal P}_{2\ell+1}^\ell(1,2\ell+1)$ is empty (there is no $(2\ell+2)$nd vertex in $P^\ell_{2\ell+1}$).

It follows that in this case the combinatorial argument expressing $\gamma(P_{2\ell+1}^\ell,x)$ in terms of $\gamma(P_m^\ell,x)$'s (with $m < 2\ell+1$) yields the relation
$$
\gamma(P_{2\ell+1}^\ell,x)-x = x\sum_{j=1}^{2\ell} \gamma(P_{2\ell+1-j}^\ell,x).
$$ 
Noting that $\gamma(P_0^\ell,x)=1$ we obtain from this that $\gamma(P_{2\ell+1}^\ell,x) = x\sum_{j=1}^{2\ell+1} \gamma(P_{2\ell+1-j}^\ell,x)$, which is \eqref{path-power-recurrence}. 

\subsubsection{The recurrence for powers of cycles} \label{subsec-proofs-power-recurrences-cycle}

Now we turn to \eqref{cycle-power-recurrence}, the recurrence relation for cycles. The proof is very similar to the proof of \eqref{path-power-recurrence}, so we aim to be brief here.

Label the vertices of $C_n$ (the cycle underlying $C_n^\ell$) $v_1, \ldots, v_n$ with $v_i$ adjacent to $v_{i+1}$ for all $i < n$, and $v_n$ adjacent to $v_1$. All dominating sets of $C_n^\ell$ have at least two vertices (a single vertex has only $2\ell$ neighbours; recall $n \geq 2\ell+2$ here), so we may partition the set of dominating sets of $C_n^\ell$ as $\cup_{f, s} {\mathcal C}_n^\ell(f,s)$, where ${\mathcal C}_n^\ell(f,s)$ is the set of dominating sets of $C_n^\ell$ in which the smallest vertex (in the natural ordering of the vertices) is $v_f$ and the second smallest is $v_{f+s}$. 

Labelling the vertices of $C_{n-s}$ (the cycle underlying $C_{n-s}^\ell$) $v^s_1, \ldots, v^s_{n-s}$ we may also partition the set of dominating sets of $C_n^\ell$ as $\cup_f {\mathcal C}_{n-s}^\ell(f)$, where ${\mathcal C}_{n-s}^\ell(f)$ is the set of dominating sets of $C_{n-s}^\ell$ in which the smallest vertex is $v^s_f$.

Using essentially the same proof as we used in the verification of \eqref{path-power-recurrence}, we get that there is a bijection from ${\mathcal C}_n^\ell(f,s)$ to ${\mathcal C}_{n-s}^\ell(f)$ that reduces the sizes of dominating sets by one. This immediately implies \eqref{cycle-power-recurrence}. The bijection is (informally) obtained by deleting all the vertices of $C_n^\ell$ that lie between $v_f$ and $v_s$, and then identifying $v_f$ and $v_s$.

\subsubsection{The recurrence for powers of paths, $n \leq 2\ell+1$} \label{subsec-proofs-power-recurrences-path-3}

We now turn to the first recurrence for powers of paths presented in Theorem \ref{thm-small-n-power-recurrences}. {\bf A5} has already been established, and {\bf A1}, {\bf A2} are easy initial conditions. For {\bf A3} and {\bf A4} we consider what happens for smaller $n$ ($2 \leq n \leq 2\ell$) in the combinatorial argument that was presented in Section \ref{subsec-proofs-power-recurrences-path-3}. As observed earlier the argument only makes sense for the set of dominating sets of $P_n^\ell$ of size at least two, and so it yields 
\begin{equation} \label{path-power-recurrence-small-n}
\gamma(P_n^\ell,x) - \gamma_1(P_n^\ell)x = x\sum_{j=1}^{n-1} \gamma(P_{n-j}^\ell,x)
\end{equation}
where recall that $\gamma_1(G)$ is the number of dominating sets of size $1$ in $G$ (note when $G$ has at least one vertex, it has no dominating sets of size $0$). 

We can easily compute $\gamma_1(P_n^\ell)$ directly:
$$
\gamma_1(P_n^\ell) = \left\{
\begin{array}{ll}
n & \mbox{if $2 \leq n \leq \ell+1$}\\
2\ell+2-n & \mbox{if $\ell+2 \leq n \leq 2\ell$}.
\end{array}
\right.
$$
Inserting this into \eqref{path-power-recurrence-small-n} we get {\bf A3} and {\bf A4} of Theorem \ref{thm-small-n-power-recurrences}.

We now turn to the second recurrence of Theorem \ref{thm-small-n-power-recurrences}. {\bf B1} is a simple initial condition, and {\bf B2} follows from the fact that for $\ell \geq n-1$, i.e. $1 \leq n \leq \ell+1$, $P_n^\ell$ is a complete graph, in which every non-empty set of vertices is dominating.

For $n \geq 2\ell+2$ (the regime of {\bf B4}) we obtain from \eqref{path-power-recurrence} that
$$
\gamma(P_{n-1}^{\ell},x)=x \sum_{i=2}^{2\ell+2}\gamma(P_{n-i}^{\ell},x).
$$ 
Rearranging and using \eqref{path-power-recurrence} we have 
\begin{eqnarray*}
\gamma(P_{n-1}^{\ell},x)-x\gamma(P_{n-2(\ell+1)}^{\ell},x) & = & x\sum_{i=2}^{2\ell+1}\gamma(P_{n-i}^{\ell},x) \\
& = & \gamma(P_n^\ell,x) - x\gamma(P_{n-1}^\ell,x).
\end{eqnarray*}
Rearranging terms again we have  
\[
\gamma(P_n^{\ell},x)=(1+x)\gamma(P_{n-1}^{\ell},x)-x \gamma(P_{n-2(\ell+1)}^{\ell},x)
\] 
for $n\geq 2\ell+2$. 
This is {\bf B4}.

We now consider what happens when 
$\ell+2 \leq n \leq 2\ell+1$ (the regime of {\bf B3}).
Since $\ell \geq 1$ we may assume $n \geq 3$. For this regime of $n$ we present a combinatorial verification of {\bf B3}. Notice that what we want to show is that for $m\geq 2$ 
\begin{equation} \label{rec1}
\gamma_m(P_{n}^{\ell})=\gamma_m(P_{n-1}^{\ell})+\gamma_{m-1}(P_{n-1}^{\ell})
\end{equation}
and also that 
\begin{equation} \label{rec2}
\gamma_1(P_{n}^{\ell})=\gamma_1(P_{n-1}^{\ell})-1.
\end{equation}
In \eqref{rec2} we are using that for $n\geq 3$ there is no empty dominating set of $P_n^{\ell}$; this fact also deals with the constant terms in {\bf B3}.

As before, let the path $P_n$ underlying $P_n^{\ell}$ have vertices $v_1,v_2,\ldots, v_n$, with $v_i$ adjacent to $v_{i+1}$ for $i=1,\ldots,n-1$. 

First note that for $m\geq n-\ell$, we have $\gamma_m(P_{n}^{\ell},x) = \binom{n}{m}$. Indeed, any set $D=\{v_{i_1}, \ldots, v_{i_m}\}$ of $m$ vertices in $P_n^\ell$ (with $i_1 < \cdots < i_m$) partitions the complement of $D$ into $m+1$ blocks --- the first block being $\{v_j: j <i_1\}$, the second block being $\{v_j: i_1 < j < i_2\}$, and so on. The sum of the sizes of the blocks is $n-m$, so each block has size at most $n-m$, which is at most $\ell$. It follows that every vertex in the first block is adjacent to $v_{i_1}$, every vertex in the second block is adjacent to both $v_{i_1}$ and $v_{i_2}$, and so on, and so $D$ is a dominating set. By Pascal's identity this yields $\gamma_m(P_{n}^{\ell})=\gamma_m(P_{n-1}^{\ell})+\gamma_{m-1}(P_{n-1}^{\ell})$ for $m\geq n-\ell$. 

What remains is to establish \eqref{rec1} for $2\leq m<n-\ell$, and also to establish \eqref{rec2}. We start with \eqref{rec1}. For each $2\leq m< n-\ell$ we partition the set of all dominating sets of $P_n^{\ell}$ of size $m$ into two classes:
     \begin{enumerate}
         \item the set ${\mathcal A}(n,\ell,m)$ of all dominating sets of $P_n^{\ell}$ of size $m$ that include $v_{n}$, and
         \item the set ${\mathcal B}(n,\ell,m)$ of all dominating sets of $P_n^{\ell}$ of size $m$ that do not include $v_{n}$. 
     \end{enumerate}
With $k=n-m-\ell$, we will show
\begin{equation} \label{u_n-in}
|{\mathcal A}(n,\ell,m)| = \gamma_{m-1}(P_{n-1}^{\ell}) + \binom{m+k-2}{m-1} 
\end{equation}
and
\begin{equation} \label{u_n-out}
|{\mathcal B}(n,\ell,m)| = \gamma_{m}(P_{n-1}^{\ell}) - \binom{m+k-2}{m-1}. 
\end{equation}
Summing \eqref{u_n-in} and \eqref{u_n-out} gives \eqref{rec1}.
     
We first prove \eqref{u_n-in}. Note that for any dominating set $D=\{v_{i_1},v_{i_2},\ldots,v_{i_{m-1}}\}$ of size $m-1$ of $P_{n-1}^{\ell}$, $D\cup \{v_n\}$ dominates $P_{n}^{\ell}$. Let $\mathcal{D}$ be the set of dominating sets $D$ of $P_n^{\ell}$ of size $m$ with $v_{n} \in D$ and with $D\setminus\{v_n\}$ not a dominating set of $P_{n-1}^{\ell}$. To establish \eqref{u_n-in} we need to show that $|\mathcal{D}|=\binom{m+k-2}{m-1}$. 

Let $D$ be in $\mathcal{D}$. Since $v_n$ is in $D$, there must exists a $D'\subseteq D$ such that $D'$ is a dominating set of $P_{n-\ell-1}^{\ell}$. Notice that for such a $D'$, $D'\cup \{v_i\}$ is a dominating set of $P_{n-1}^{\ell}$ for any $v_i$ in $\{v_{n-\ell},\ldots, v_{n-1}\}$. So, using that $D\setminus\{v_n\}$ doesn't dominate $P_{n-1}^{\ell}$, we see that all the $m-1$ vertices in $D\setminus\{v_n\}$ must come from $\{v_1,\ldots, v_{\ell+1}\}$. Moreover, none of the vertices $v_{i}$ with $m+k-1\leq i\leq n-1$ can be in $D$ (recall $k=n-m-\ell$). This says that $D\setminus\{v_n\}\subseteq\{v_1,\ldots, v_{m+k-2}\}$ (here the final index really should be $\min\{m+k-2,\ell+1\}$; but note since $n \leq 2\ell+1$ and $m+k=n-\ell$, we have $m+k-2\leq \ell+1$). 

We have shown $|\mathcal{D}| \leq \binom{m+k-2}{m-1}$. To see the equality, note that any $m-1$ vertices in $\{v_1,\ldots, v_{m+k-2}\}$ will dominate $P_{n-\ell-1}^{\ell}$ but won't dominate $P_{n-1}^{\ell}$ (because $k\leq \ell$). This concludes the proof of \eqref{u_n-in}.

We now move on to \eqref{u_n-out}. Again, we see that every dominating set of $P_{n}^{\ell}$ that does not include $v_n$ is a dominating set of $P_{n-1}^{\ell}$. So we want to count the number of dominating sets $D$ of size $m$ of $P_{n-1}^{\ell}$ that do not dominate $P_{n}^{\ell}$. A necessary condition here is that $D\subseteq \{v_1,\ldots,v_{n-1-j}\}$, and more precisely we must have $D\subseteq\{v_1,\ldots, v_{m+k-1}\}$ and $v_{m+k-1}\in D$ (with again $k=n-m-j$). As in the verification of \eqref{u_n-out}, since any $m-1$ vertices of $\{v_1,\ldots, v_{m+k-2}\}$ dominate $P_{m+k-1}^{\ell}$ (because $j\geq k$) we get that the number of dominating sets $D$ of size $m$ of $P_{n-1}^{\ell}$ that don't dominate $P_{n}^{\ell}$ is equal to $\binom{m+k-2}{m-1}$. This establishes \eqref{u_n-out}, and so \eqref{rec1}.

We now move on to \eqref{rec2}. All dominating sets of $P_{n}^{\ell}$ of size $1$ are dominating sets of $P_{n-1}^{\ell}$. The only dominating set of $P_{n-1}^{\ell}$ of size $1$ that doesn't dominate $P_{n}^{\ell}$ is $\{v_{n-1-\ell}\}$. This yields \eqref{rec2}, and so completes the verification of {\bf B3}.

\subsubsection{A shorter proof of \eqref{path-power-recurrence} when $n \geq 3\ell+2$} \label{subsec-proofs-power-recurrences-path-2} 

We note, as an aside, that there is a particularly clean way to present the above argument that establishes \eqref{path-power-recurrence} for $n \geq 3\ell+2$. Say that $D$ is a {\it relaxed dominating set} in $P_n^\ell$ if every vertex in $P_n^\ell$, other than possibly the first $\ell$ vertices of the underlying path $P_n$, is either in $D$ or adjacent to something in $D$. Denote by  $\gamma^r(P_n^\ell,k)$ the number of relaxed dominating sets of size $k$, and let
$$
\gamma^r(P_n^\ell,x) = \sum_{k \geq 0} \gamma^r(P_n^\ell,k)x^k.
$$

We can express $\gamma(P_n^\ell,x)$ in terms of the $\gamma^r(P_m^\ell,x)$'s. Note that any dominating set of $P_n^\ell$ must include at least one of the first $\ell+1$ vertices of the underlying path $P_n^\ell$ (to ensure that the first vertex of the path is either in the set or adjacent to something in the set). Let the $\ell+1$ initial vertices of the path be $v_1, \ldots, v_{\ell+1}$. For each $i=1, \ldots, \ell+1$ there is a one-to-one correspondence between on the one hand 
\begin{quote}
dominating sets of $P_n^\ell$ in which the first vertex from among $v_1, \ldots, v_{\ell+1}$ that is in the set is $v_i$
\end{quote}
and on the other hand
\begin{quote}
relaxed dominating sets of $P_{n-i}^\ell$. 
\end{quote}
The correspondence is  
obtained by deleting the first $i$ vertices of $P_n^\ell$ to go from $P_n^\ell$ to $P_{n-i}^\ell$, and removing $v_i$ from the dominating set of $P_n^\ell$ to get a relaxed dominating set of $P_{n-i}^\ell$. This correspondence drops the size of a set by one, and so it follows that   
\begin{equation} \label{relaxed1}
\gamma(P_n^\ell,x) = x\sum_{i=1}^{\ell+1} \gamma^r(P_{n-i}^\ell,x).
\end{equation}
Note that \eqref{relaxed1} makes sense for $n \geq \ell+1$. 

Similarly we can express $\gamma^r(P_n^\ell,x)$ in terms of the $\gamma^r(P_m^\ell,x)$'s. Note that any relaxed dominating set of $P_n^\ell$ must include at least one of the first $2\ell+1$ vertices of the underlying path $P_n^\ell$ (to ensure that the $(\ell+1)$st vertex of the path is either in the set or adjacent to something in the set). Let the $2\ell+1$ initial vertices of the path be $v_1, \ldots, v_{2\ell+1}$. For each $i=1, \ldots, 2\ell+1$ there is a one-to-one correspondence between on the one hand 
\begin{quote}
relaxed dominating sets of $P_n^\ell$ in which the first vertex from among $v_1, \ldots, v_{2\ell+1}$ that is in the set is $v_i$ 
\end{quote}
and on the other hand 
\begin{quote}
relaxed dominating sets of $P_{n-i}^\ell$. 
\end{quote}
The correspondence is obtained by deleting the first $i$ vertices of $P_n^\ell$ to go from $P_n^\ell$ to $P_{n-i}^\ell$, and removing $v_i$ from the relaxed dominating set of $P_n^\ell$ to get a relaxed dominating set of $P_{n-i}^\ell$. This correspondence drops the size of a set by one, and so it follows that   
\begin{equation} \label{relaxed2}
\gamma^r(P_n^\ell,x) = x\sum_{i=1}^{2\ell+1} \gamma^r(P_{n-i}^\ell,x).
\end{equation}
Note that \eqref{relaxed2} makes sense for $n \geq 2\ell+1$. 

Using \eqref{relaxed1} for the first equation below, \eqref{relaxed2} for the second, then changing the order of summation, and finally using \eqref{relaxed1} in the fourth equation, we obtain
\begin{eqnarray*}
\gamma(P_n^\ell,x) & = & x\sum_{i=1}^{\ell+1} \gamma^r(P_{n-i}^\ell,x) \\
& = & x^2\sum_{i=1}^{\ell+1} \sum_{j=1}^{2\ell+1} \gamma^r(P_{n-i-j}^\ell,x) \\
& = & x^2\sum_{j=1}^{2\ell+1} \sum_{i=1}^{\ell+1} \gamma^r(P_{n-i-j}^\ell,x) \\
& = & x\sum_{j=1}^{2\ell+1} \gamma(P_{n-j}^\ell,x).
\end{eqnarray*}
Note that this is valid for $n \geq 3\ell+2$.

\subsection{Unimodality of sequences of polynomials} \label{subsec-proofs-bb-gen}

We now turn to the proof of Theorem \ref{thm-bb-gen}. For each fixed $k \geq 3$ we proceed by induction on $n$, with $n=k$ as the base case. 

For the inductive step, let us assume that ${\mathcal P}_n$ holds for some $n \geq k$. Our goal is to establish ${\mathcal P}_{n+1}$. 

That the coefficients of $f_{n+1}$ are non-negative follows immediately from the recurrence relation satisfied by the $f_i(x)$'s.  

For $0 \leq r \leq n+1$ and $c \geq 0$ let $a_{r,c}$ be the coefficient of $x^c$ in $f_r(x)$. By the induction hypothesis each of the sequences $(a_{r,c})_{c \geq 0}$ is non-negative and unimodal for $0 \leq r \leq n$, and there is a barely increasing sequence $(m_r)_{r=0}^n$ such that $a_{r,m_r}$ is a mode of the sequence $(a_{r,c})_{c \geq 0}$.

That each $a_{n+1,c}$, for $c \geq 0$, is non-negative follows immediately from the recurrence relation satisfied by the $f_i(x)$'s. Thus to establish ${\mathcal P}_{n+1}$, we need to show that $(a_{n+1,c})_{c \geq 0}$ is unimodal with a mode at $c=m_n$ or $m_n+1$.  

We begin by observing that $(a_{n+1,c})_{c \geq m_n+1}$ is monotone decreasing. Indeed, for $t \geq m_n+1$ we have
\begin{eqnarray*}
a_{n+1,t} & = & \sum_{j=1}^k a_{n-k+1,t-1} \\
& \geq & \sum_{j=1}^k a_{n-k+1,t} \\
& = & a_{n+1,t+1}.
\end{eqnarray*}
The equalities above are simply applications of the recurrence relation satisfied by the $f_i(x)$'s. The inequality follows from the fact that $(a_{n-k+1,c})_{c \geq 0}$ has a mode at $m_{n-k+1}$, and by the increasing property of the sequence of modes we have $m_{n-k+1} \leq m_n \leq t-1$, so $a_{n-k+1,t-1} \geq a_{n-k+1,t}$. 

It thus remains to show that $(a_{n+1,c})_{c \leq m_n}$ is monotone increasing. Let $0 \leq j \leq k-1$ be such that $m_{n-k+1}=m_n-j$ (note $j \geq 0$ because the sequence of modes is increasing, and $j \leq k-1$ because the sequence of modes is barely increasing). By an almost identical argument to the one used to show that $(a_{n+1,c})_{c \geq m_n+1}$ is monotone decreasing, we see that $(a_{n+1,c})_{c \leq m_n-j+1}$ is monotone increasing.

To complete the argument that $(a_{n+1,c})_{c \leq m_n}$ is monotone increasing, fix $i$ with $m_n-j+1 \leq i \leq m_n-1$. We need to show that $a_{n+1,i} \leq a_{n+1,i+1}$, which is equivalent to
\begin{equation} \label{int1}
\sum_{p=1}^k a_{n+1-p,i} \geq \sum_{p=1}^k a_{n+1-p,i-1}.  
\end{equation}
Note that because $i \leq m_n-1$ we have $a_{n,i+1} \geq a_{n,i}$. By the recurrence this says that
\begin{equation} \label{int2}
\left(\sum_{p=2}^k a_{n+1-p,i}\right) + a_{n-k,i} \geq \left(\sum_{p=1}^k a_{n+1-p,i-1}\right) + a_{n-k,i-1}. 
\end{equation}
Now because $i \geq m_n-j+1$, so $i-1 \geq m_n-j$, and by increasing property of the sequence of modes (which says that $m_{n-k}\leq m_n-j$) we have that $a_{n-k,i} \leq a_{n-k,i-1}$. Inserting into \eqref{int2} we get
\begin{equation} \label{int3}
\left(\sum_{p=2}^k a_{n+1-p,i}\right) \geq \left(\sum_{p=1}^k a_{n+1-p,i-1}\right).
\end{equation}
Finally noting that $a_{n,i} \geq a_{n,i-1}$ (since $i \leq m_n-1$) we deduce \eqref{int1} from \eqref{int3}, completing the induction. 

\subsection{Unimodality for powers of paths and cycles} \label{subsec-proofs-powers-of-path-cycle}

Here we give the proof of Theorem \ref{thm-powers-of-path-cycle} --- unimodality of the domination sequence of powers of paths and cycles --- beginning with the case of paths. 

We start with large $\ell$, specifically $\ell \geq 9$. Note that $\gamma(P_0^\ell,x)=1$, which has a unimodal coefficient sequence with a mode at $\lceil 0/2 \rceil$, while for $1 \leq n \leq \ell+1$ we have that $P_n^\ell$ is simply the complete graph on $n$ vertices, and so
$$
\gamma(P_n^\ell,x) = (1+x)^n - 1,
$$
which for each $n$ has unimodal coefficient sequence with a mode at $\lceil n/2 \rceil$.

For $n \geq \ell+2$ the graph $P_n^\ell$ has minimum degree $\ell$. As long as $n \leq 2^{\ell/2}$ we get from Theorem \ref{thm-from-bb-min-degree} that $\gamma(P_n^\ell,x)$
has unimodal coefficient sequence with a mode at $\lceil n/2 \rceil$. Noting that $(\lceil n/2 \rceil)_{n \geq 0}$ is barely increasing, and bearing the recurrence \eqref{path-power-recurrence} in mind, we conclude that as long as $\lfloor 2^{\ell/2}\rfloor \geq 2\ell+1$ we can apply Theorem \ref{thm-bb-gen} to conclude that the coefficient sequence of $\gamma(P_n^\ell,x)$ is unimodal for all $n \geq 0$. This establishes the path case of Theorem \ref{thm-powers-of-path-cycle} for all $\ell \geq 9$.

To complete the path case of Theorem \ref{thm-powers-of-path-cycle} we simply have to establish that for each $1 \leq \ell \leq 8$, each of the polynomials in the sequence $(\gamma(P_n^\ell,x))_{n=0}^{2\ell+1}$ has a unimodal coefficient sequence, and that the family has a sequence of modes that is barely increasing. Using either of the recurrences given in Theorem \ref{thm-small-n-power-recurrences} this finite task is easily accomplished (we used {\tt Mathematica} for this computation), and specifically it is easily established that each polynomial under examination has a mode at $\lceil n/2 \rceil$.  

We now turn to the cycle case of Theorem \ref{thm-powers-of-path-cycle}. By \eqref{eq-small-n-cycle-recurrence} we have that for $1 \leq n \leq 2\ell+1$ each of $\gamma(C_n^\ell,x)$ has non-negative and unimodal coefficient sequence, and has a mode at $\lceil n/2 \rceil$. By \eqref{cycle-power-recurrence} and \eqref{eq-small-n-cycle-recurrence} we have
\begin{eqnarray}
\gamma(C_{2\ell+2}^\ell,x) & = &\sum_{i=1}^{2\ell+1} \left((1+x)^i-1\right) \nonumber \\
& = & (1+x)^{2\ell+2} - (2\ell+2)x - 1. \label{cycle_2l+2}
\end{eqnarray}
(We can also see \eqref{cycle_2l+2} directly --- all subsets of the vertex set of $C_{2\ell+2}^\ell$ of size two or more are dominating sets, but no subset of size $1$ or $0$ is.) So $\gamma(C_{2\ell+2}^\ell,x)$ also has non-negative and unimodal coefficient sequence with a mode at $\ell+1$. By Theorem \ref{thm-bb-gen} we conclude that $\gamma(C_n^\ell,x)$ is unimodal for all $n$.

\end{document}